\documentclass[11 pt]{amsart}

\usepackage{times}
\usepackage{geometry}
\usepackage{amssymb}
\usepackage{latexsym,amssymb,  amsmath, amscd, amsfonts}
\usepackage{graphicx}
\usepackage[percent]{overpic}
\usepackage{pdfsync}
\usepackage{units}
\usepackage{hyperref}
\usepackage{euscript}
\usepackage{multicol}
\usepackage{epstopdf}
\usepackage{paralist}

\newtheorem{theorem}{Theorem}

\newtheorem{proposition}[theorem]{Proposition}

\theoremstyle{definition}
\newtheorem{definition}[theorem]{Definition}

\newtheorem{remark}[theorem]{Remark}

\def\fig#1{Figure~\ref{fig:#1}}
\def\prop#1{Proposition~\ref{prop:#1}}

\newcommand{\R}{\mathbb{R}}      %From Becky's file

\def\lk{{\operatorname{Lk}}}

\def\rib{{\operatorname{Rib}}}
\def\cr{{\operatorname{Cr}}}
\def\len{{\operatorname{Len}}}

\begin{document}

\title[Folded Ribbon Knots]{Folded ribbon knots in the plane}

\author[Denne]{Elizabeth Denne}
\address{Elizabeth Denne, Washington \& Lee University, Department of Mathematics, Lexington VA}
\date{June 28, 2018.}                                           % Activate to display a given date or no date

\begin{abstract}
This survey reviews Kauffman's model of folded ribbon knots: knots made of a thin strip of paper folded flat in the plane. The ribbonlength is the length to width ratio of such a ribbon, and the ribbonlength problem asks to minimize the ribbonlength for a given knot type. We give a summary of known results. For the most part, these are upper bounds of ribbonlength of twist knots and certain families of torus knots. We discuss result of G.~Tian \cite{Tian}, which give upper bounds of ribbonlength in terms of crossing number. In addition, it turns out the choice of fold affects the ribbonlength. We end with a discussion of three different types of folded ribbon equivalence and give examples illustrating their relationship to ribbonlength.
\end{abstract}

\maketitle

\section{Introduction}\label{section:intro}

\begin{center}
\begin{figure}[htbp]
\includegraphics{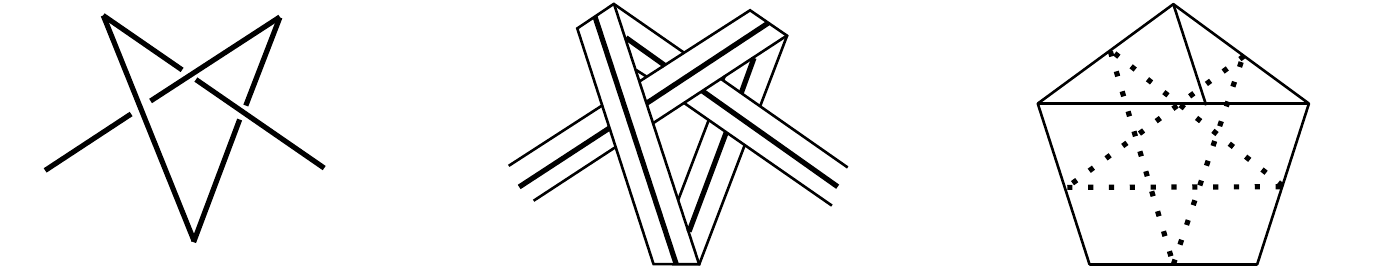}
\caption{Creating a folded ribbon trefoil knot.}
\label{fig:trefoil-ribbon}
\end{figure}
\end{center}

We can create a ribbon knot in $\mathbb R^3$ by taking a long, rectangular piece of paper, tying a knot in it, and connecting the two ends. We then flatten the ribbon into the plane, origami style, with folds in the ribbon appearing only at the corners.  In \fig{trefoil-ribbon} left and center, we show part of a trefoil knot and the corresponding folded ribbon trefoil. If we join the two ends of the ribbon and tighten, we get  a ``tight'' folded ribbon trefoil shown on the right.

 Such a {\em folded ribbon knot} was first modeled by L. Kauffman \cite{Kauf05}. (He called them flat knotted ribbons.) Kauffman defined the {\em ribbonlength} (the length to width ratio) of folded ribbons knots, and asked to find the least length of ribbon needed for a knot type, given a choice of width. This is known as the {\em (folded) ribbonlength problem}. A minimal ribbonlength folded ribbon knot can be considered to be folded ``tightly''. This idea is neatly illustrated with the construction of a tight folded ribbon trefoil knot in \fig{trefoil-ribbon} right, which has a pentagon as its boundary shape. 
Kauffman~\cite{Kauf05} gave the conjectured minimal ribbonlength for this and a tight figure eight knot.   

Understanding folded ribbonlength reveals interesting relationships between geometry and topology, and there are natural connections between folded ribbons and other areas of mathematics and science. For example, the ribbonlength problem may also be thought of as a 2-dimensional analogue of the {\em ropelength problem}: that of finding the minimum amount of rope needed to tie a knot in a rope of unit diameter. (See for instance \cite{BS99,CKS,gm,lsdr,Pie}.)  Folded ribbon knots arise naturally from considering (smooth) ribbon knots in space. These are used, for example, to model cyclic duplex DNA in molecular biology with the two boundaries corresponding to the two edges of the DNA ladder (see for instance \cite{Adams,BF}).  Folded ribbon knots have connections to other parts of knot theory as well. We will see later that grid diagrams of knots, and mosaic knots, can easily produce folded ribbon knots. Grid diagrams of knots have been extensively studied \cite{OSZ} and are used, for example, in the combinatorial formulation of knot Floer homology. Mosaic knots  \cite{LK,KS} are used to model quantum knots  which describe a physical quantum system. 

In the past, there have been a number of recreational articles about tying knots in strips of paper. In particular, \cite{Adams, Ashley, CR, John, Wel} all described the construction of a pentagon as the boundary shape of a folded ribbon trefoil knot. Some of these also found other regular $n$-gons; for example D.A. Johnson~\cite{John} gave the construction of a regular hexagon by folding a trivial 2-component link in a certain way.  Constructing any regular polygons as the boundary of a folded ribbon knot appeared to be a harder problem. This was finally solved in the 1999 master's thesis of  L. DeMaranville~\cite{DeM}. She described which $(p,q)$ torus knots\footnote{Here, a $(p,q)$ torus knot is assumed to wrap $p$ times around the meridional direction and $q$ times around the longitude direction of a torus.} can be easily converted to a folded ribbon knot, and showed how to build all regular $n$-gons for $n>6$, by tying certain families of folded torus ribbon knots. While there is often more than one way to construct a particular $n$-gon, the $(p,2)$ torus knots  give all odd $p$-gons, and the $(q+1,q)$ torus knots give all $2q$-gons (here $p, 2q>6$). An interesting open question is to find what other shapes can be formed by torus links, and what can be said about their folded ribbonlengths.  

%%%%%%%%%%%%%%%%%%%%%%%%%%%
%%%%%%%%%%%%%%%%%%%%%%%%%%%%%%

\section{Modeling Folded Ribbon Knots}\label{section:model}

Before we proceed further, let us begin by reviewing some familiar definitions (see for instance  \cite{Adams, Kauf-book, Liv-book}). A {\em tame knot} is an embedding of $S^1$ in $\R^3$ which is ambient isotopic to a polygonal knot. A {\em link} is a disjoint union of knots, and we will use the word knot to mean either a knot or a link. A {\em projection of a knot K} is the image of {\em K} under a projection from $\mathbb R^3$ to a plane, and a {\em knot diagram} adds gaps in a knot projection to show over- and under-crossing information.

\begin{center}
\begin{figure}[htbp]
\includegraphics{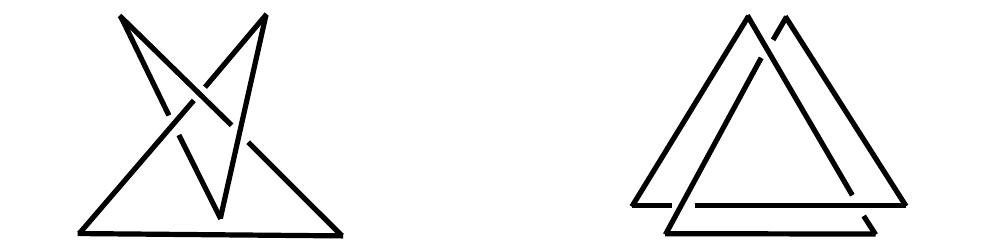}
\caption{Polygonal knot diagrams of the trefoil knot with five and six edges.}
\label{fig:trefoil-diagrams}
\end{figure}
\end{center}

A formal definition of a folded ribbon knot can be found in \cite{Rib-Smith, DKTZ}. In this section, we provide enough details to give the reader the big picture. When we model a folded ribbon knot, we view the knot as a {\em polygonal knot diagram}.  \fig{trefoil-diagrams} shows two different polygonal knot diagrams for the trefoil knot. 
A polygonal knot diagram has a finite number of {\em vertices} denoted by $v_1, ..., v_n$, and {\em edges}~$e_i$ defined by  $e_1=[v_1, v_2]$, \dots , $e_n=[v_n,v_1]$. If the knot diagram is oriented, then we assume that the labeling follows the orientation.

We pause to note that we do not require our polygonal knot diagrams to be {\em regular}. For example, take a ribbon which is an annulus and then fold it flat with just two folds.  The polygonal knot diagram is made of two edges and we understand that one edge always lies over the other so that the crossing information is consistent. We will refer to this diagram as the 2-stick unknot. Recall that the {\em stick index} of a knot $K$ is defined to be the least number of line segments needed to construct a polygonal embedding of $K$ in $\R^3$ (see \cite{Adams}). We can define the {\em projection stick index} to be the minimum number of sticks needed for a polygonal knot diagram of $K$. We have just seen the unknot has projection stick index two, and regular stick index three. 
Together with undergraduate students, C. Adams~\cite{Adams-Shayler, Adams-PS} showed that the projection stick index of the trefoil knot is five, while the regular stick index is six (illustrated in \fig{trefoil-diagrams}). Indeed, we expect the projection stick index to be smaller than the stick index since the edges in the knot diagram are not rigid sticks in space, they have crossing information instead. 

\begin{figure}[htbp]
\begin{center}
\begin{overpic}{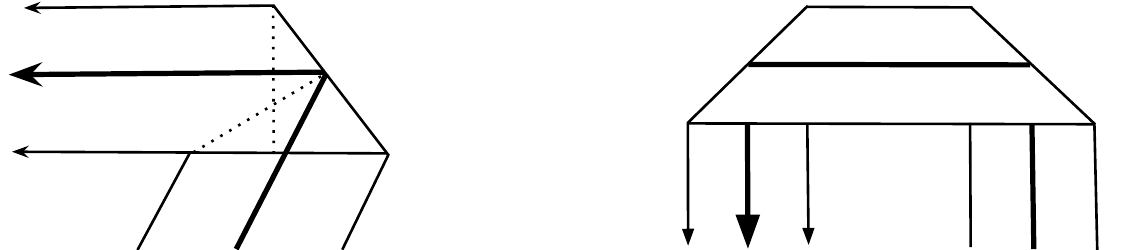}
\put(23,23){$A$}
\put(35,8){$B$}
\put(30,15.5){$C=v_i$}
\put(24.8,13){$\theta_i$}
\put(26,6){$E$}
\put(12.5,6){$D$}
\put(21,16.5){$F$}
\put(21,6){$G$}
\put(5,17){$e_i$}
\put(15,1){$e_{i-1}$}

\put(80,18){$e_i$}
\put(92,-2){$e_{i-1}$}
\put(67,-1){$e_{i+1}$}
\put(92,17){$v_i$}
\put(61,18){$v_{i+1}$}
\end{overpic}
\caption{On the left, a close-up view of a ribbon fold. On the right, the construction of the ribbon centered on edge $e_i$. }
\label{fig:ribbon-construct}
\end{center}
\end{figure}
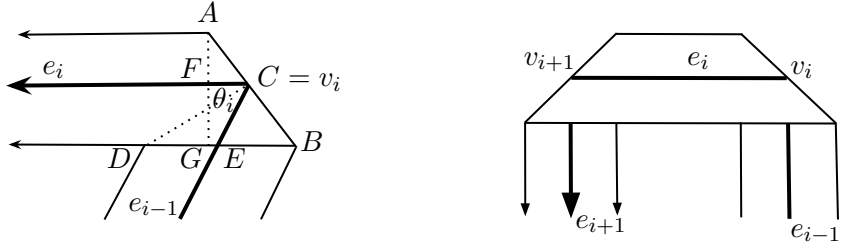

To construct a width $w$ folded ribbon knot, we view a polygonal knot diagram $K$ as the centerline of the ribbon. Then  the fold lines of the ribbon are perpendicular to the angle bisectors at each of the knot diagram's vertices.  

\begin{definition} Given an oriented polygonal knot diagram $K$, we define the {\em fold angle} at vertex $v_i$ to be the angle $\theta_i$ (where $0\leq\theta_i\leq \pi$) between edges $e_{i-1}$ and $e_i$.  Then, the {\em oriented folded ribbon knot of width} $w$, denoted $K_w$, is constructed as follows:
 \begin{enumerate}

\item  First, construct the fold lines. At each vertex $v_i$ of $K$, find the fold angle~$\theta_i$. If $\theta_i<\pi$, place a fold line of length $w/\cos(\frac{\theta_i}{2})$ centered at $v_i$ perpendicular to the angle bisector of $\theta_i$.   If $\theta_i=\pi$, there is no fold line.
\item Second, add in the ribbon boundaries. For each edge $e_i$, join the ends of the fold lines at  $v_i$ and $v_{i+1}$. By construction, each boundary line is parallel to, and distance $w/2$ from $K$.
\item The ribbon inherits an orientation from $K$.
\end{enumerate}
\end{definition}

This construction is illustrated in \fig{ribbon-construct}. On the left, the fold angle is $\theta_i=\angle ECF$, the angle bisector is $DC$, and the fold line is $AB$.  Using the geometry of the figure we see that $\angle GAB=\theta_i/2$ in right triangle $\triangle AGB$. Thus $|AB|=w/\cos(\frac{\theta_i}{2})$ guarantees the ribbon width $|AG|=w$.  We say that  the fold angle is $\theta_i=\angle ECF$ is {\em positive}, since $e_i$ is to the left of $e_{i-1}$. If $e_i$ were to the right of $e_{i-1}$, then it would be {\em negative}.

\begin{figure}[htbp]
\begin{center}
\begin{overpic}{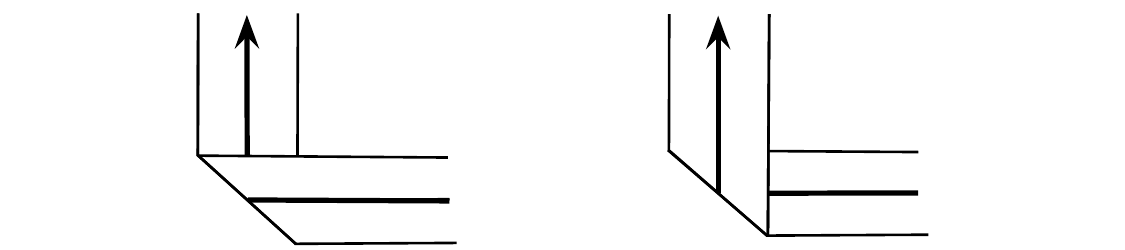}
\put(36, 6){$e_{i-1}$}
\put(23, 20){$e_{i}$}
\put(17, 4){$v_i$}

\put(76, 6.5){$e_{i-1}$}
\put(65, 20){$e_{i}$}
\put(60, 4){$v_i$}
\end{overpic}
\caption{A right underfold (left) and a right overfold (right).}
\label{fig:folding-info}
\end{center}
\end{figure}
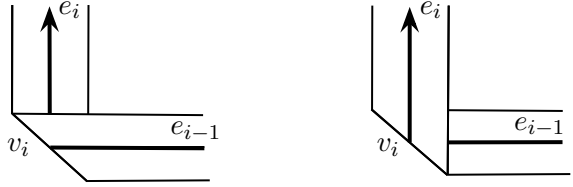

Observe that near a fold line, there is a choice of which ribbon lies above the other. Thus a polygonal knot diagram with $n$ vertices has $2^n$ possible folded ribbon knots depending on the choices made.  There is an {\em overfold} at vertex $v_i$ if the ribbon corresponding to segment $e_{i}$ is over the ribbon of segment $e_{i-1}$ (see \fig{folding-info} right). Similarly, there is an {\em underfold} if the ribbon corresponding to $e_{i}$ is under the ribbon of $e_{i-1}$. 

\begin{definition} The choice of overfold or underfold at each vertex of $K_w$ is called the {\em folding information}, and is denoted by $F$.
\end{definition}

In our construction of a folded ribbon knot,  we start with a polygonal knot diagram then build the ribbon. There appears to be no restriction placed on the width, and yet very wide ribbons might not be physically possible. For example in \fig{ribbon-construct} (right), the ribbon width can not be more than the length of $e_i$. To guarantee physically possible ribbons, we require that the ribbon has consistent crossing information.

\begin{definition} Given an oriented knot diagram $K$, we say the folded ribbon $K_{w,F}$ of width $w$ and folding information $F$ is {\em allowed} provided
\begin{enumerate}
\item The ribbon has no singularities (is immersed), except at the fold lines.
\item $K_w$ has consistent crossing information, and moreover this agrees
\begin{enumerate} \item with the folding information given by $F$, and
\item with the crossing information of the knot diagram $K$.
\end{enumerate}
\end{enumerate}
\end{definition}

When a folded ribbon $K_{w,F}$ is allowed, the consistent crossing information means that a straight ribbon segment cannot ``pierce'' a fold, it either lies entirely above or below the fold, or lies between the two ribbons segments joined at the fold.  From now on {\bf we assume that our folded ribbon knots have an allowed width}. This is a reasonable assumption, since we can always construct folded ribbon knots for ``small enough'' widths. 

 \begin{proposition}[\cite{Rib-Smith}] Given any regular polygonal knot diagram $K$ and folding  information $F$, there is a constant $C>0$ such that an allowed folded ribbon knot $K_{w,F}$ exists for all $w<C$. \qed
\end{proposition}

%%%%%%%%%%%%%%%%%%%%%%%%%%%%%%%%%%%

\section{Ribbonlength} \label{section:ribbonlength}

Given a particular folded ribbon knot, it is very natural to wonder what is the least length of ribbon needed to tie it.  More formally, we define a scale invariant quantity, called {\em ribbonlength}, as follows.

\begin{definition} [\cite{Kauf05}] The {\em (folded) ribbonlength}, $\rib(K_{w,F})$, of a folded ribbon knot $K_{w,F}$ is the quotient of the length of $K$ to the width $w$:
$$\rib(K_{w,F})=\frac{\len(K)}{w}.$$
\end{definition}     

When minimizing the ribbonlength of a folded ribbon knot, we have two choices. We can fix the width and minimize the length, or, we can can fix the length and maximize the width. As mentioned above, Kauffman \cite{Kauf05} gave the conjectured minimal ribbonlength for the trefoil and figure eight knots.  In 2008, B. Kennedy, T.W. Mattman, R. Raya and D. Tating~\cite{KMRT} gave upper bounds on the ribbonlength of the $(p,2)$,  $(q+1,q)$,  and $(2q+1,q)$ families of torus knots, 
using ideas in the Master's theses of DeMaranville~\cite{DeM} and Kennedy~\cite{Ken}.
They did not expect the bounds to be minimal, and in fact gave shorter versions of the $(5,2)$ and $(7,2)$ torus knots.  

An interesting open question is to understand the relationship between the ribbonlength of a knot $K$ and its crossing number\footnote{The crossing number of a knot is the minimum number of crossings in any regular knot diagram of the knot.} $\cr(K)$. In particular, to find constants $c_1, c_2, \alpha, \beta$ such that 
\begin{equation}\label{bounds}
c_1\cdot( \cr(K))^\alpha\leq \rib(K_w)\leq c_2\cdot (\cr(K))^\beta.
\end{equation}
R. Kusner \cite{Kauf05} conjectured that ribbonlength has upper and lower bounds that are linear in the crossing number, that is $\alpha=\beta=1$ in Equation~\ref{bounds}.
Kennedy {\it et al.}~\cite{KMRT} made a first pass at the bounds in Equation~\ref{bounds} by using the fact the crossing number of a $(p,q)$ torus knot is $\min\{p(q-1),q(p-1)\}$. This allowed them to show the ribbonlength's upper bound is quadratic in the crossing number for the $(p,2)$ torus knots, and linear for the $(q+1,q)$ and $(2q+1,q)$ torus knots.

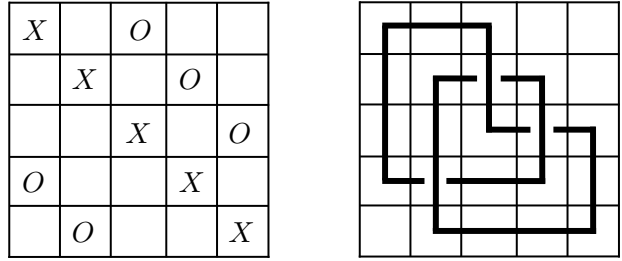
\begin{figure}[htbp]
\begin{center}
\begin{overpic}{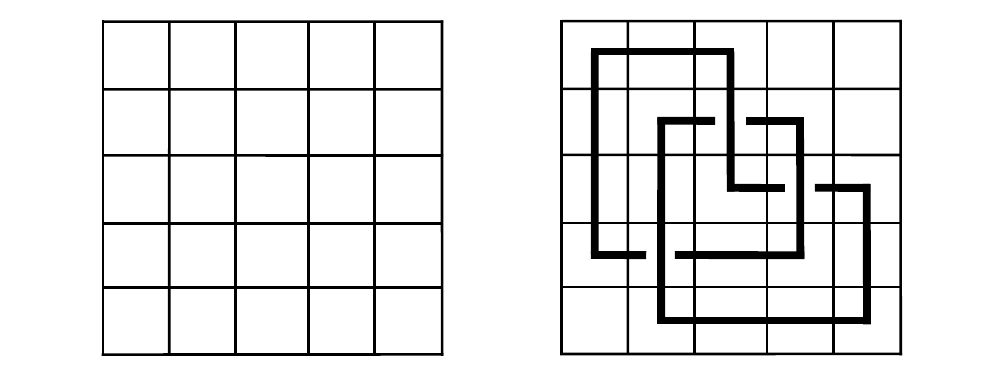}
\put(12,30.5){$X$}
\put(26,30.5){$O$}
\put(18.5,24){$X$}
\put(32.5,24){$O$}
\put(25.5,17){$X$}
\put(39,17){$O$}
\put(32.5,10.5){$X$}
\put(12,10.5){$O$}
\put(39,4){$X$}
\put(18.5,4){$O$}
\end{overpic}
\end{center}
\caption{A grid diagram (left) and corresponding knot diagram (right) of a trefoil knot.}
\label{grid-diagram}
\end{figure}

In 2017, G. Tian~\cite{Tian} used grid diagrams of knots to make further progress. A {\em grid diagram}, with grid number $n$, is an $n\times n$ square grid with $n$ $X$'s and $n$ $O$'s arranged so that every row and column contains exactly one $X$ and one $O$. A grid diagram of the trefoil knot is given on the left in Figure~\ref{grid-diagram}. It turns out that every knot has a grid diagram associated to it \cite{Cro, OSZ}. 

Now, any grid diagram gives a knot diagram in a standard way: connect $O$ to $X$ in each row, connect $X$ to $O$ in each column, and have the vertical line segments always cross over the horizontal ones. This process is illustrated on the right in Figure~\ref{grid-diagram}.
From here, we can see that any grid diagram of a knot can be used to create a folded ribbon knot whose width is the sidelength of the squares in the grid.

\begin{figure}[htbp]
\begin{center}
\begin{overpic}[scale=0.8]{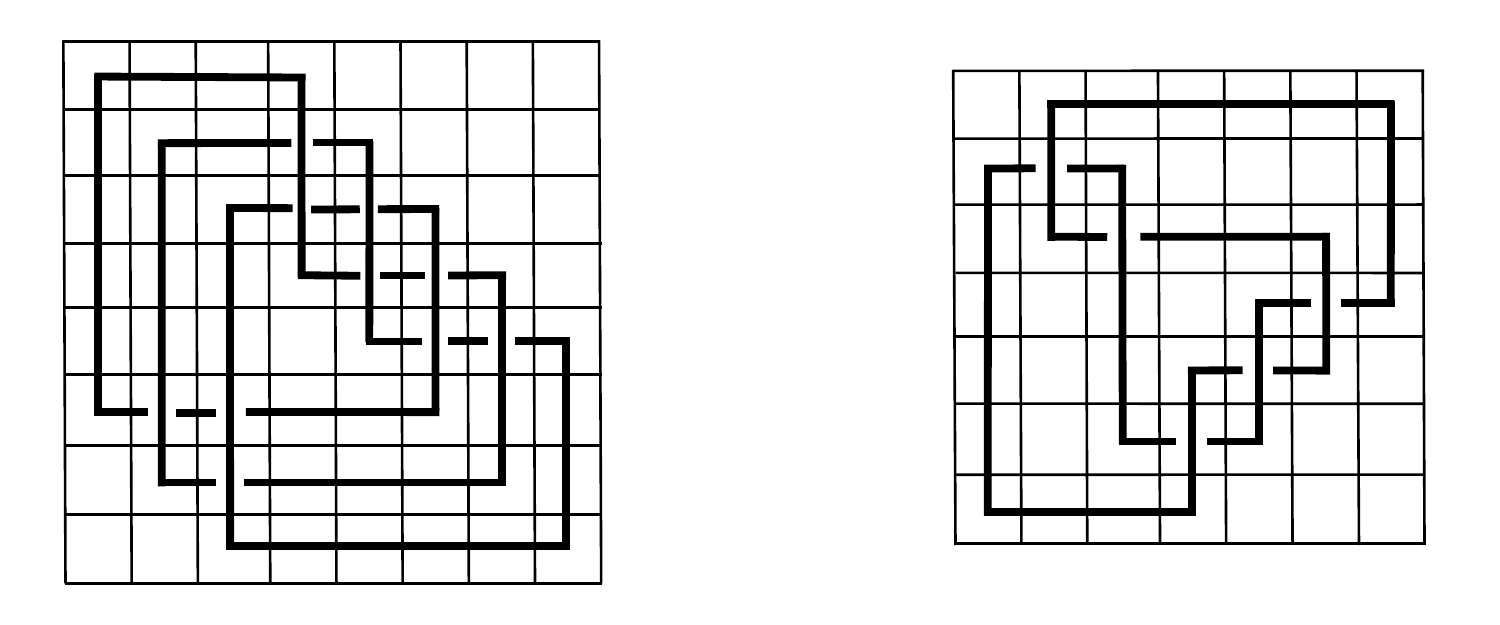}
\end{overpic}
\end{center}
\caption{The grid diagram of the $(5,3)$ torus knot on the left, and the $3$-twist knot on the right.}
\label{grid-examples}
\end{figure}

The {\em grid index}, $g(K)$, is the minimal grid number, and others have proved \cite{BP, HKON} that the grid index is bounded above by crossing number: $g(K)\leq \cr(K)+2$.   Tian~\cite{Tian} argued that given a knot $K$ represented by a $g(K)\times g(K)$ grid diagram, we can estimate the length of the folded ribbon knot. Each horizontal distance between the $X$ and $O$ is at most $g(K)-1$, hence the sum of all horizontal distance is at most $g(K)(g(K)-1)$. The same is true for vertical distances. Thus we obtain
$$\rib(K)\leq 2g(K)(g(K)-1)\leq 2(\cr(K)+1)(\cr(K)+2)\leq 12(\cr(K))^2.
$$
Tian then used a grid diagram of $n$-twist knots $T_n$, to show the ribbonlength of the corresponding folded ribbon $T_n$ knot is $2(4n+8)$. Since the crossing number of $T_n$ is $n+2$, we find $\rib(T_n)\leq 8\cr(T_n)$.
A similar argument for $(p,q)$ torus knots shows the ribbonlength  is also bounded above by $8\cr(K)$.   Figure~\ref{grid-examples} shows the kind of grid diagrams for torus and twist knots used to obtain the bounds.

Altogether, we know that $\beta\leq 2$ in Equation~\ref{bounds}, and we have several families of knots where $\beta=1$. Tian's results both improves and extends the results of Kennedy {\em et al.}~\cite{KMRT} described above.
So far, no one has proved any lower bounds on ribbonlength in terms of crossing number. We end this section by noting that progress has been made on finding upper and lower bounds of ropelength in terms of crossing number (see for instance \cite{BS99, BS07, CFM, CKS, DEKZ, DEPZ, DEY}).

%%%%%%%%%%%%%%%%%%%%%%%%%%
\section{Ribbon equivalence and ribbonlength}\label{section:equivalence}

As we saw in Section~\ref{section:model}, there is a choice of folding information at each vertex of the knot diagram. It turns out that the folding information affects the ribbonlength.  To get a handle on these differences, we give three different definitions of ribbon equivalence (see also \cite{Rib-Smith, DKTZ}). We first begin by defining ribbon linking number.

\begin{figure}[htbp]
\begin{center}
\begin{overpic}{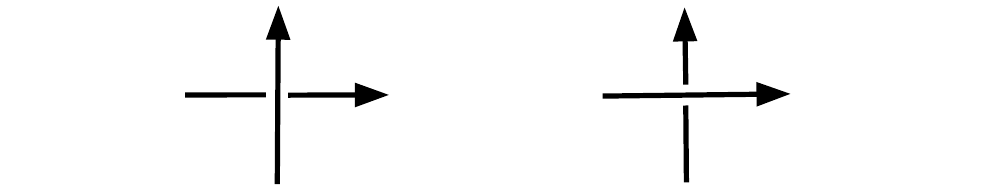}
\end{overpic}
\end{center}
\caption{The crossing on the left is labelled -1, the crossing on the right +1.}
\label{crossingvalue}
\end{figure}

The {\em linking number} is an invariant from knot theory (see for instance \cite{Adams, Kauf-book, Liv-book}) used to determine the degree to which components of a link are joined together.  Given an oriented two component link $L=A\cup B$, the linking number  $\lk(A,B)$ is defined to be one half the sum of $+1$ crossings and $-1$ crossings between $A$ and $B$. (See Figure~\ref{crossingvalue}.)

Although we have described the construction of folded ribbon knots in $\mathbb R^2$, we can also consider the ribbons that these diagrams represent in $\mathbb R^3$. That is, as {\em framed knots}. The ribbon linking number was defined for these ribbons  (see \cite{Cal59, Cal61, Kauf-book}), but it equally applies to our situation.

 \begin{definition}
 Given an oriented folded ribbon knot $K_{w,F}$, we define the {\em (folded) ribbon linking number} to be the linking number between the knot diagram and one boundary component of the ribbon.  We denote this as $\lk(K_{w,F})$, or $\lk(K_w)$.
 \end{definition}

We are now ready to define three different kinds of ribbon equivalence, starting with the most restrictive. 

\begin{definition}(Link equivalence)
Two oriented folded ribbon knots are {\em (ribbon) link equivalent} if they have equivalent knot diagrams with the same ribbon linking number.
\end{definition}

\begin{figure}[htbp]
\begin{center}
\begin{overpic}{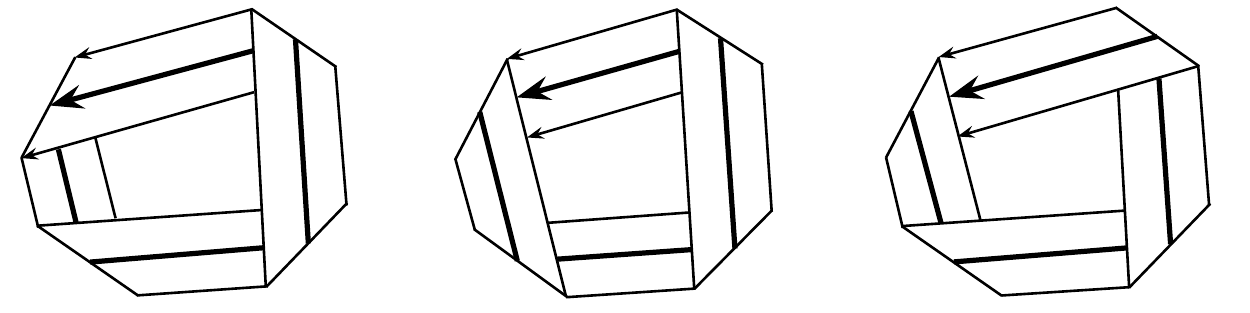}
\end{overpic}
 \caption{The left and center 4-stick folded ribbon unknots have ribbon linking number 0, while the one on the right has ribbon linking number $-2$. } 
\label{fig:4unknot}
\end{center}
\end{figure}

For example, the left and center folded ribbon unknots in \fig{4unknot} are link equivalent (with ribbon linking number 0), while the one on the right is not link equivalent to them (with ribbon linking number $-2$). This example shows that there can be different looking folded ribbon knots with the same ribbon linking number.

\begin{definition}(Topological equivalence)  Two oriented  folded ribbon knots are {\em topologically (ribbon) equivalent} if they have equivalent knot diagrams and, when considered as ribbons in $\mathbb R^3$, both ribbons are topologically equivalent to a M\"obius strip or both ribbons are topologically equivalent to an annulus.  
\end{definition}

For example, all of the 4-stick folded ribbon unknots in \fig{4unknot} are topologically equivalent.

\begin{definition}(Knot diagram equivalence) Two folded ribbon knots are {\em knot diagram equivalent} if they have equivalent knot diagrams.
\end{definition}

\begin{figure}[htbp]
\begin{center}
\includegraphics{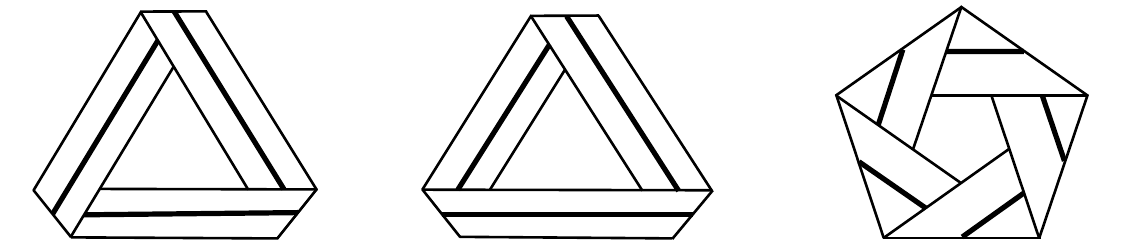}
\caption{On the left and the center, the two different kinds of folding information for 3-stick unknots. On the right, a 5-stick folded ribbon unknot. }
\label{fig:unknots-ribbonlength}
\end{center}
\end{figure}

For example, the 3-stick and 4-stick folded ribbon unknots in Figures~\ref{fig:4unknot} and \ref{fig:unknots-ribbonlength}  are knot diagram equivalent, but are not topologically equivalent, nor link equivalent. 
 
 \begin{remark} The {\em ribbonlength problem} asks us to minimize the ribbonlength of a folded ribbon knot, while staying in a fixed topological knot type. That is, with respect to knot diagram equivalence of folded ribbon knots.  We can also ask to minimize the ribbonlength of folded ribbon knots with respect to topological and link equivalence.
\end{remark}

%%%%%%%%%

We pause to remark that the previous work on ribbonlength \cite{Kauf05,KMRT,Tian} found upper bounds on the ribbonlength with respect to knot diagram equivalence. Together with undergraduate students \cite{Rib-Smith, DKTZ}, we have made a first pass at finding bounds on ribbonlength with respect to topological and link equivalence. We start by considering unknots.

Any polygonal unknot diagram can be reduced to a 2-stick unknot, and the width of such a diagram can be made large as we like. Thus the minimum ribbonlength of any unknot with respect to knot diagram equivalence is 0.  If we minimize ribbonlength with respect to topological equivalence, then we have already considered the topological annulus in the 2-stick unknot. The 3-stick unknot is a M\"obius band, and this is straightforward to understand.

\begin{proposition}[\cite{DKTZ} ] The minimum ribbonlength of an 3-stick folded ribbon unknot $K_{w,F}$ is less than or equal to
\begin{enumerate}
\item $3\sqrt{3}$ \ when the folds are all the same type (\fig{unknots-ribbonlength} left),
\item $\sqrt{3}$ \ when one fold is of different type to the other two (\fig{unknots-ribbonlength} center). \qed
\end{enumerate}
\end{proposition} 

It turns out we can do even better and show that the equilateral triangle gives the ribbonlength minimizer in the first case. 

\begin{theorem}[\cite{DKTZ}]
The minimum ribbonlength for the 3-stick folded unknot is $3\sqrt{3}$ where all folds have the same folding information, and occurs when the knot diagram is an equilateral triangle. \qed
\end{theorem}

Thus the ribbonlength of an unknot is less than or equal to $\sqrt{3}$ when the ribbon is topologically equivalent to a M\"obius strip, and is 0 when the ribbon is topologically equivalent to an annulus. When ribbon linking number is taken into consideration the situation is more complex. For example, for 3-stick unknots with ribbon linking number $\pm 3$ the minimum ribbonlength is  equal to $3\sqrt{3}$, and is less than or equal to $\sqrt{3}$ for 3-stick unknots with ribbon linking number $\pm 1$. What about higher linking numbers? Assume that all folds of the an $n$-stick folded ribbon unknot are the same. Then for odd $n\geq 4$ the ribbon linking number is $n$, while for even $n\geq 4$ the ribbon linking number is $n/2$.  In either case we can get an upper bound on the ribbonlength of $n$-stick unknots by considering the case where the knot diagrams are regular $n$-gons (as in \fig{unknots-ribbonlength} right). 

\begin{proposition}[\cite{DKTZ}] \label{prop:min-unknot} The ribbonlength of an $n$-stick folded ribbon unknot (for $n\ge4$) is less than or equal to $n \cot(\frac{\pi}{n})$. \qed
\end{proposition} 

Thus \prop{min-unknot} gives a reasonable upper bound on ribbonlength with respect to link equivalence for $n$-stick unknots with $n$ odd sides and ribbon linking number $n$. Just how far can we improve our ribbonlength bounds (with respect for link equivalence) for the unknots with any ribbon linking number?  

%%%%%%%%

\begin{figure}[htbp]
\begin{center}
\includegraphics[scale=0.8]{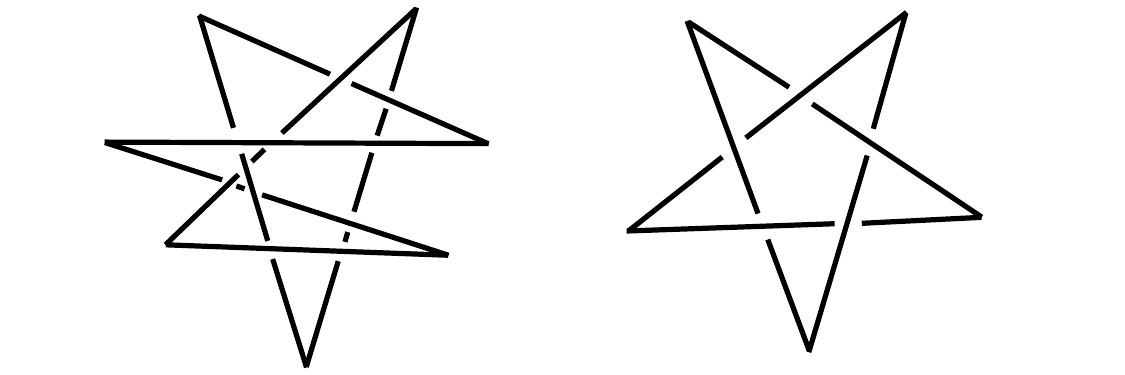}
\caption{Two different polygonal diagrams of the $(5,2)$ torus knot. The one on the left has a smaller ribbonlength than the one on the right, despite having 2 more edges.}
\label{fig:T52-equivalent}
\end{center}
\end{figure}

Moving now from unknots to nontrivial knots, it is natural to wonder about the relationship between the ribbonlength of a knot, the number of edges in the knot diagram, and the projection stick index of the knot. Recall that Kennedy {\em et al.}~\cite {KMRT}, found that there was a smaller ribbonlength for the $(5,2)$ and $(7,2)$ torus knots, simply by adding two more edges to the knot diagram and rearranging. The knot diagram for the $(5,2)$ torus knot with smaller ribbonlength is shown on the left in \fig{T52-equivalent}. In \cite{DKTZ}, we showed that the two folded ribbon $(5,2)$ torus knots corresponding to the knot diagrams in  \fig{T52-equivalent} are not ribbon link equivalent. We found a sequence of Reidemeister moves connecting the two knot diagrams, and we showed the corresponding folded ribbons differed by a full twist (regardless of the starting folding information). This example shows two interesting things. Firstly, the minimal ribbonlength of a knot (with respect to knot diagram equivalence) does not necessarily occur when the knot diagram has the projection stick index. Secondly, the candidate for minimal ribbonlength is {\em not} ribbon link equivalent to the folded ribbon knot whose diagram has the projection stick index. 
There is much to think about here!

In summary, we have seen that folded ribbon knots are both interesting to study in their own right, and have many connections to other parts of knot theory. We close by listing just some of the open questions that have come up in this survey of folded ribbon knots.

\begin{itemize}
\item Prove that the ribbonlength of the trefoil knot is minimized by the configuration with pentagonal boundary given in \fig{trefoil-ribbon} (right).
\item Improve the upper and lower bounds relating ribbonlength to crossing number given in Equation~\ref{bounds}.
\item Is a minimal ribbonlength folded ribbon knot always given by a knot diagram that has more edges than the projection stick index?
\item Is a minimal ribbonlength folded ribbon knot ever ribbon link equivalent to a folded ribbon knot whose diagram has the projection stick index?
\item What are the ribbon linking numbers generated by a knot diagram with the projection stick index? How is this related to minimum ribbonlength?
\item Minimize ribbonlength  of nontrivial knots with respect to  link equivalence.
\end{itemize}

%%%%%%%%%%%%%%%%%%%%%%%%

\section{Acknowledgments}
I wish to thank Jason Cantarella for many helpful math conversation, and John M. Sullivan and Nancy Wrinkle for discussions on smooth ribbon knots immersed in the plane \cite{DSW}. 

I have worked with undergraduate students over many years on the folded ribbonlength problem. Special thanks go to former Smith College students: Shivani Aryal, Eleanor Conley,  Shorena Kalandarishvili, Emily Meehan and Rebecca Terry. Special thanks also go to former Washington \& Lee Students:  Mary Kamp, and Catherine (Xichen) Zhu.

%%%%%%%%%%%%%%%%%%%%%%%%%%

\bibliographystyle{amsalpha}

\end{document}